\documentclass [11pt,twoside,a4paper]{article}
\usepackage{amsfonts}
\usepackage{amsthm}
\usepackage{amsmath}
\usepackage{amstext}
\usepackage{amssymb}
\usepackage{mathrsfs}
\usepackage{amscd}
\usepackage{xypic}
\usepackage{epsf}              
\usepackage{graphicx}          
\usepackage{fancybox}          
\usepackage{color}             
\usepackage{fancyhdr}
\usepackage[hang,footnotesize]{caption2}  

\def\mathcal{\mathscr}
\newfont{\aaa}{cmb10 at 19pt}
\newfont{\bbb}{cmb10 at 11pt}

\def\v1{\vspace{1mm}}

\def\leq{\leqslant}

\def\geq{\geqslant}

\newcommand{\beq}{\begin{equation}}
\newcommand{\eeq}{\end{equation}}
\newcommand{\bey}{\begin{eqnarray}}
\newcommand{\eey}{\end{eqnarray}}
\newcommand{\beyy}{\begin{eqnarray*}}
\newcommand{\eeyy}{\end{eqnarray*}}

\makeatletter
\def\@evenhead{
\vbox{\hbox to \textwidth {}{\hspace{0mm}{\footnotesize
\thepage}}{\hspace{8cm} {\footnotesize {Manli Song, Wenjuan Li
}}} \protect\vspace{1truemm}\relax \hrule depth0pt
height0.15truemm width\textwidth}}
\def\@evenfoot{}
\def\@oddhead{\vbox{\hbox to \textwidth
{{\hspace{0cm}{\footnotesize Weighted Caffarelli-Kohn-Nirenberg type inequalities related to Grushin type operators}\hfill{\footnotesize
\thepage}}\hspace{0mm}}{} \protect\vspace{1truemm}\relax\hrule
depth0pt height0.15truemm width\textwidth}}
\def\@oddfoot{}
\makeatother

\begin{document}



\setcounter{page}{1}
\qquad\\[8mm]

\noindent{\aaa{Weighted Caffarelli-Kohn-Nirenberg type inequalities related to Grushin type operators}}\\[1mm]

\noindent{\bbb Manli Song,\quad Wenjuan Li}\\[-1mm]

\noindent\footnotesize{School of Natural and Applied
Sciences, Northwestern Polytechnical University, Xi'an, Shaanxi 710129, China}\\[6mm]

\normalsize\noindent{\bbb Abstract}\quad We consider the Grushin type operator on $\mathbb{R}^{d}_x \times \mathbb{R}^{k}_y$ with the form
\begin{equation*}
G_\mu=\overset{d}{\underset{i=1}{\sum}}\partial_{x_i}^2+\left(\overset{d}{\underset{i=1}{\sum}}x_i^2\right)^{2\mu}\overset{k}{\underset{j=1}{\sum}}\partial_{y_j}^2.
\end{equation*}
and derive weighted Hardy-Sobolev type inequalities and weighted Caffarelli-Kohn-Nirenberg type inequalities related to $G_\mu$.\vspace{0.3cm}

\footnotetext{
Corresponding author: Wenjuan Li, E-mail:
liwj@nwpu.edu.cn}

\noindent{\bbb Keywords}\quad Grushin type operator, Weighted Hardy-Sobolev inequality, Weighted Caffarelli-Nirenberg type inequality.\\
{\bbb MSC}\quad 26D10, 35H20\\[0.4cm]

\newtheorem{theorem}{Theorem}[section]
\newtheorem{preliminaries}{Preliminaries}[section]
\newtheorem{definition}{Difinition}[section]
\newtheorem{main result}{Main Result}[section]
\newtheorem{lemma}{Lemma}[section]
\newtheorem{proposition}{Proposition}[section]
\newtheorem{corollary}{Corollary}[section]
\newtheorem{remark}{Remark}[section]
\section{Introduction}
Hardy-Sobolev inqualities and Caffarelli-Kohn-Nirenberg inequalities on the Euclidean space play an important role in mathematics and applied fields. They are very useful tools to study various interesting problems in partial differential equations, such as eigenvalue problems, existence problems of equation with singular weights, regularity problems, etc.

The initial work of first order interpolation inequalities with weights (Caffarelli-Kohn-Nirenberg inequality) was given by Caffarelli et al \cite{CKN}. The result is stated as follows:
\begin{theorem} Let $p$, $q$, $r$, $\alpha$, $\beta$, $\gamma$, $\sigma$ and $a$ satisfy
\begin{equation}\label{real-condition}
\left\{
\begin{aligned}
&p,q\geq 1, r>0, 0\leq a\leq 1\\
&\frac{1}{p}+\frac{\alpha}{n}, \frac{1}{q}+\frac{\beta}{n}, \frac{1}{r}+\frac{\gamma}{n}>0,
\end{aligned}
\right.
\end{equation}
where $\gamma=a\sigma+(1-a)\beta$. Then there exists a positive constant $C$ such that the following inequality holds for all $u\in C_0^\infty(\mathbb{R}^n)$
\begin{equation}\label{CKN-inequality}
\Big\||x|^\gamma u\Big\|_{L^r}\leq C\Big\||x|^\alpha|\nabla u|\Big\|^a_{L^p}\left\||x|^\beta u\right\|_{L^q}^{1-a}
\end{equation}
if and only if the following relations hold:
\begin{equation}\label{dimensional-balance}
\frac{1}{r}+\frac{\gamma}{n}=a(\frac{1}{p}+\frac{\alpha-1}{n})+(1-a)(\frac{1}{q}+\frac{\beta}{n})
\end{equation}
(this is dimensional balance),
\begin{align*}
&0\leq \alpha-\sigma\quad \text{if} \quad a>0\\
&\alpha-\sigma\leq1\quad \text{if} \quad a>0\quad\text{and}\quad \frac{1}{p}+\frac{\alpha-1}{n}=\frac{1}{r}+\frac{\gamma}{n}.
\end{align*}
Furthermore, on any compact set in parameter space in which \eqref{real-condition}, \eqref{dimensional-balance} and $0\leq\alpha-\sigma\leq1$ hold, the constant $C$ is bounded.
\end{theorem}

When $a=1$, we see that \eqref{CKN-inequality} are reduced to Hardy-Sobolev inequalities, i.e., Hardy-Sobolev inequalities are special cases of Caffarelli-Kohn-Nirenberg inequalities. Later, Lin \cite{Lin} generalized \eqref{CKN-inequality} to cases including derivatives of any order. Badiale and Tarantello \cite{BT} derived a class of more general Hardy-Sobolev inequalities with singular weights depending only on partial variables. Recently, Hardy-Sobolev type inequalities have been extended to noncommutative field vectors. In the Heisenberg group setting, we refer the readers to see \cite{DAM}\cite{Han}\cite{HN}\cite{HNZ}\cite{HZD}, ect.

In this paper, we shall prove weighted Hardy-Sobolev type inequalities and weighted Caffarelli-Kohn-Nirenberg type inequalities related to the Grushin type operator \begin{equation*}
G_\mu=\Delta_x+|x|^{2\mu}\Delta_y,
\end{equation*}
where $x\in\mathbb{R}^d$ and $y\in\mathbb{R}^k$. Let $Q=d+(1+\mu)k$ be the homogeneous dimension and
 \begin{equation*}
\rho=\rho(x,y)=\Big(|x|^{2+2\mu}+(1+\mu)^2|y|^2\Big)^{\frac{1}{2+2\mu}}
\end{equation*}
be the distance function from the origin $(x,y)$ on $\mathbb{R}^{d}_x \times \mathbb{R}^{k}_y$ to and $\nabla_\mu=(\nabla_x, |x|^\mu\nabla_y)$ be the gradient operator. We have $G_\mu=\nabla_\mu\cdot\nabla_\mu$.

For the Grushin type operator, D' Ambrosio \cite{DAM1} proved some Hardy type inequalities and gave sharp estimates in some cases. Here, we recall a result in \cite{DAM1}:
\begin{theorem} \label{Hardy-type-theorem}
Let $p>1$ and $\alpha\in\mathbb{R}$ satisfy $\frac{1}{p}+\frac{\alpha}{Q}>0$.  Then there exists a positive constant $C=\left(\frac{p}{Q-p+\alpha p}\right)^p$ such that for any $u\in D_\alpha^{1,p}(\mathbb{R}^{d+k})$,
\begin{equation}\label{Hardy-type}
\int_{\mathbb{R}^{d+k}}\rho^{\alpha p}\frac{|x|^{\mu p}}{\rho^{\mu p}}\frac{|u|^p}{\rho^p}dxdy\leq C\int_{\mathbb{R}^{d+k}}\rho^{\alpha p}|\nabla_\mu u|^pdxdy,
\end{equation}
where $D_\alpha^{1,p}(\mathbb{R}^{d+k})$ denotes the closure of $C_0^\infty(\mathbb{R}^{d+k})$ with respect to the norm
\begin{equation*}
\|u\|_\alpha^{1,p}=\left(\int_{\mathbb{R}^{d+k}}|\rho^\alpha\nabla_\mu u|^pdxdy\right)^{\frac{1}{p}}.
\end{equation*}
\end{theorem}

 Niu and Dou \cite{ND} established Hardy-Sobolev inequalities related to $G_\mu$. Zhang et al \cite{ZHD} obtained a class of weighted Hardy-Sobolev inequalities and a class of weighted Caffarelli-Kohn-Nirenberg inequalities in the special case $p=2$. The weighted Hardy-Sobolev inequalities are listed as follows:
\begin{theorem}(see \cite{ZHD})
If $0\leq s\leq2<Q$, $\alpha>\frac{2-Q}{2}$, there exists a positive constant $C=C(s, \alpha, \mu, Q)$ such that for any  $u\in D_\alpha^{1,2}(\mathbb{R}^{d+k})$,
\begin{equation*}
\int_{\mathbb{R}^{d+k}}\frac{|x|^{\mu s}}{\rho^{\mu s}}\frac{|\rho^\alpha u|^{2_*(s,2,Q)}}{\rho^s}dxdy\leq C\left(\int_{\mathbb{R}^{d+k}}|\rho^\alpha\nabla_\mu u|^2dxdy\right)^{\frac{Q-s}{Q-2}},
\end{equation*}
where we put $p_*(s,p,Q)=\frac{p(Q-s)}{Q-p}$ for any $1<p<Q$.
\end{theorem}

The weighted Caffarelli-Kohn-Nirenberg inequalities related to $G_\mu$ for the case $p=2$ are given:
\begin{theorem}(see \cite{ZHD}) Let $q$, $r$, $\alpha$, $\beta$, $\gamma$, $\sigma$ and $a$ satisfy
\begin{equation*}
\left\{
\begin{aligned}
&q\geq 1, r>0, 0\leq a\leq 1\\
&d+\mu(\alpha-\gamma)r>0, d+\mu(\alpha-\beta)q>0\\
&\gamma r+Q>0, \beta q+Q>0, 2\alpha+Q>0,
\end{aligned}
\right.
\end{equation*}
where $\gamma=a\sigma+(1-a)\beta$. Then there exists a positive constant $C$ such that the following inequality holds for all $u\in C_0^\infty(\mathbb{R}^{d+k})$
\begin{align*}
&\left(\int_{\mathbb{R}^{d+k}}\left(\frac{|x|^\mu}{\rho^\mu}\right)^{(\alpha-\gamma)r}\rho^{\gamma r}|u|^rdxdy\right)^{\frac{1}{r}} \\
\leq &C
\left(\int_{\mathbb{R}^{d+k}}\rho^{2\alpha}|\nabla_\mu u|^2dxdy\right)^{\frac{a}{2}}
\left(\int_{\mathbb{R}^{d+k}}\left(\frac{|x|^\mu}{\rho^\mu}\right)^{(\alpha-\beta)q}\rho^{\beta q}|u|^qdxdy\right)^{\frac{1-a}{q}}
\end{align*}
if and only if the following relations hold:
\begin{equation*}
\frac{1}{r}+\frac{\gamma}{Q}=a(\frac{1}{2}+\frac{\alpha-1}{Q})+(1-a)(\frac{1}{q}+\frac{\beta}{Q})
\end{equation*}
(this is dimensional balance),
\begin{align*}
&0\leq \alpha-\sigma\quad \text{if} \quad a>0\\
&\alpha-\sigma\leq1\quad \text{if} \quad a>0\quad\text{and}\quad \frac{1}{2}+\frac{\alpha-1}{Q}=\frac{1}{r}+\frac{\gamma}{Q}.
\end{align*}
\end{theorem}

The paper is to establish the following weighted Caffarelli-Kohn-Nirenberg inequalities related to $G_\mu$ for $1<p<Q$:
\begin{theorem} \label{CKN-Grushin} Let $p$, $q$, $r$, $\alpha$, $\beta$, $\gamma$, $\sigma$ and $a$ satisfy
\begin{equation*}
\left\{
\begin{aligned}
&1<p<Q, q\geq 1, r>0, 0\leq a\leq 1\\
&d+\mu(\alpha-\gamma)r>0, d+\mu(\alpha-\beta)q>0\\
&\alpha p+Q>0, \beta q+Q>0, \gamma r+Q>0,
\end{aligned}
\right.
\end{equation*}
where
\begin{equation}\label{alpha-beta-gamma}
\gamma=a\sigma+(1-a)\beta.
\end{equation}
Then there exists a positive constant $C$ such that the following inequality holds for all $u\in C_0^\infty(\mathbb{R}^{d+k})$
\begin{align}\label{CKN-inequality-new}
&\left(\int_{\mathbb{R}^{d+k}}\left(\frac{|x|^\mu}{\rho^\mu}\right)^{(\alpha-\gamma)r}\rho^{\gamma r}|u|^rdxdy\right)^{\frac{1}{r}} \nonumber\\
\leq &C
\left(\int_{\mathbb{R}^{d+k}}\rho^{\alpha p}|\nabla_\mu u|^p dxdy\right)^{\frac{a}{p}}
\left(\int_{\mathbb{R}^{d+k}}\left(\frac{|x|^\mu}{\rho^\mu}\right)^{(\alpha-\beta)q}\rho^{\beta q}|u|^qdxdy\right)^{\frac{1-a}{q}}
\end{align}
if and only if the following relations hold:
\begin{equation}\label{dimensional-balance-new}
\frac{1}{r}+\frac{\gamma}{Q}=a(\frac{1}{p}+\frac{\alpha-1}{Q})+(1-a)(\frac{1}{q}+\frac{\beta}{Q})
\end{equation}
(this is dimensional balance),
\begin{align}
&0\leq \alpha-\sigma\quad \text{if} \quad a>0 \label{index}\\
&\alpha-\sigma\leq1\quad \text{if} \quad a>0\quad\text{and}\quad \frac{1}{p}+\frac{\alpha-1}{Q}=\frac{1}{r}+\frac{\gamma}{Q}.\label{case>1}
\end{align}
\end{theorem}

To prove Theorem \ref{CKN-Grushin}, by employing the idea in \cite{Han}\cite{ZHD}, we first need to obtain a class of weighted Hardy-Sobolev type inequalities for $1<p<Q$:
\begin{theorem}\label{W-H-S-new}
If $1<p<Q$, $0\leq s\leq p$, $\alpha>\frac{p-Q}{p}$, there exists a positive constant $C=C(s,p,\alpha,\mu,Q)$ such that for any  $u\in D_\alpha^{1,p}(\mathbb{R}^{d+k})$,
\begin{equation}\label{weighted-H-S-Grushin}
\int_{\mathbb{R}^{d+k}}\frac{|x|^{\mu s}}{\rho^{\mu s}}\frac{|\rho^\alpha u|^{p_*(s,p,Q)}}{\rho^s}dxdy\leq C\left(\int_{\mathbb{R}^{d+k}}|\rho^\alpha\nabla_\mu u|^pdxdy\right)^{\frac{Q-s}{Q-p}}.
\end{equation}
\end{theorem}

\begin{remark} \label{remark}When $a=1$, conditions of Theorem \ref{CKN-Grushin} imply
\begin{equation*}
0\leq \alpha-\sigma=\alpha-\gamma\leq1,\quad \frac{1}{r}+\frac{\sigma}{Q}= \frac{1}{p}+\frac{\alpha-1}{Q},
\end{equation*}
and then $p\leq r\leq p^*=\frac{Qp}{Q-p}$. Therefore, there exists $t\in[0,1]$ satisfying
\begin{equation*}
r=tp+(1-t)p^*=\frac{p(Q-tp)}{Q-p},\quad \text{and}\quad (\alpha-\sigma)r=tp.
\end{equation*}
Replacing $r$ and $\alpha-\sigma$ into \eqref{CKN-inequality-new}, it is easy to see that \eqref{CKN-inequality-new} is reduced to \eqref{weighted-H-S-Grushin}, which is exactly a weighted Hardy-Sobolev type inequality.
\end{remark}

\begin{remark} Since the methods based on radial symmetry in \cite{CKN} are no longer suitable for $G_\mu$, Zhang et al \cite{ZHD} adopted a different idea. They first proved weighted Hardy-Sobolev type inequalities related to $G_\mu$ for the case $p=2$ and then derived the associated weighted Caffarelli-Kohn-Nirenberg inequalities for the special case. Inspired by their work, we extended their results to all the cases $1<p<Q$. However, it is still open for the cases $p\geq Q$.
\end{remark}

This paper is organized as follows. The next section introduces some definitions and basic facts related to $G_\mu$. In Section 3, we establish a Sobolev-Stein embedding theorem and Hardy-Sobolev type inequalities related to $G_\mu$. Furthermore, we prove Theorem \ref{W-H-S-new}. The final section is devoted to the proof of Theorem \ref{CKN-Grushin}.
\\[4mm]
\section{Preliminary}
We shall introduce some notions and basic facts about the Grushin type operators. Let $\mu$ be a positive real number and $(x,y)\in\mathbb{R}_x^d\times\mathbb{R}_y^k=\mathbb{R}^{d+k}$ with $d,k\geq1$. We denote by $|x|$ (resp. $|y|$) the Euclidean norm in $\mathbb{R}^d$ (resp. $\mathbb{R}^k$), i.e., $|x|^2=\overset{d}{\underset{i=1}{\sum}}x_i^2$ (resp. $|y|^2=\overset{k}{\underset{j=1}{\sum}}y_j^2$).
The symbol $\nabla_x$ (resp. $\nabla_y$) and $\Delta_x$ (resp. $\Delta_y$) stand respectively for the usual gradient operator and the Laplace operator on $\mathbb{R}^d$ (resp. $\mathbb{R}^k$).

The Grushin type vector fields are defined by
\begin{equation*}
X_i=\partial_{x_i},\,Y_j=|x|^\mu\partial_{y_j},\, i=1, 2,\cdots, d,\,j=1,2,\cdots,k,
\end{equation*}
and the corresponding gradient operator and divengent operator are denoted respectively by
\begin{align*}
&\nabla_\mu=(X_1, X_2, \cdots, X_d, Y_1, Y_2, \cdots, Y_k)=(\nabla_x, |x|^\mu\nabla_y),\\
&\text{div}_\mu(u_1,u_2,\cdots,u_{d+k})=\overset{d}{\underset{i=1}{\sum}}X_iu_i+\overset{k}{\underset{j=1}{\sum}}Y_ju_{j+d}.
\end{align*}
Denote the Grushin type operator by
\begin{equation*}
G_\mu=\overset{d}{\underset{i=1}{\sum}}X_i^2+\overset{k}{\underset{j=1}{\sum}}Y_j^2=\Delta_x+|x|^{2\mu}\Delta_y=\nabla_\mu\cdot\nabla_\mu.
\end{equation*}

A family of the dilations $\{\delta_\lambda: \lambda>0\}$ on $\mathbb{R}^{d+k}$ is defined by
\begin{equation*}
\delta_\lambda(x,y)=(\lambda x, \lambda^{1+\mu} y),\,(x,y)\in\mathbb{R}^{d+k},
\end{equation*}
and $Q=d+(1+\mu)k$ is the corresponding homogeneous dimension. It is easy to see that the vector fields $X_i$ and $Y_j$ are homogeneous of degree one with respect to the dilation, i.e., $X_i(\delta_\lambda)=\lambda \delta_\lambda(X_i)$, $Y_j(\delta_\lambda)=\lambda \delta_\lambda(Y_j)$, and hence $\nabla_\mu(\delta_\lambda)=\lambda\delta_\lambda(\nabla_\mu)$ and $G_\mu(\delta_\lambda)=\lambda^2\delta_\lambda(G_\mu)$.

The distance function from the origin to $(x,y)$ on $\mathbb{R}^{d+k}$ is defined by
\begin{equation*}
\rho=\rho(x,y)=\left(\left(\overset{d}{\underset{i=1}{\sum}}x_i^2\right)^{1+\mu}+(1+\mu)^2\overset{k}{\underset{j=1}{\sum}}y_j^2\right)^{\frac{1}{2+2\mu}}=\Big(|x|^{2+2\mu}+(1+\mu)^2|y|^2\Big)^{\frac{1}{2+2\mu}}.
\end{equation*}
It is not difficult to check that $\rho$ is homogeneous of degree one with respect to $\delta_\lambda$ and
\begin{equation}\label{gradient-Delta}
|\nabla_\mu\rho|=\frac{|x|^\mu}{\rho^\mu},\,G_\mu\rho=(Q-1)\frac{|x|^{2\mu}}{\rho^{2\mu+1}}.
\end{equation}
Furthermore, $\Gamma=C_\mu\rho^{2-Q}$ is the fundamental solution at the origin of $G_\mu$ (see \cite{DL}).

Denote the open ball of radius $R$ centered at the origin by
\begin{equation*}
B_R=\{(x,y)\in\mathbb{R}^{d+k}: \rho(x,y)<R\}.
\end{equation*}

Recalling the explicit polar transform defined by D'Ambrosio \cite{DAM}, one has
\begin{equation*}
dxdy=\rho^{Q-1}d\rho d\sigma,
\end{equation*}
where $d\sigma=\left(\frac{1}{1+\mu}\right)^k|\sin \theta|^{\frac{d}{2}-1}|\cos \theta|^{k-1}d\theta d\omega_dd\omega_k$, $\omega_d$ and $\omega_k$ denote respectively the usual surface measures on $\mathbb{R}^d$ and $\mathbb{R}^k$. In addition, the criteria for the integrability of $|x|^p\rho^q$ was given as follows:

(1) if $p+d>0$ and $p+q+Q>0$, then $\int_{B_2}|x|^p\rho^qdxdy<+\infty$;

(2) if $p+d>0$ and $p+q+Q<0$, then $\int_{\mathbb{R}^{d+k}\setminus B_1}|x|^p\rho^qdxdy<+\infty$.
\noindent\\[4mm]
\section{Proof of Theorem \ref{W-H-S-new}}
Firstly, we need to prove the Sobolev-Stein embedding result related to Grushin type operators.
\begin{theorem} \label{Sobolev-Stein-Embedding} If $1<p<Q$, there exists a positive constant $C=C(p,\mu, Q)$ such that for any  $u\in D_0^{1,p}(\mathbb{R}^{d+k})$,
\begin{equation*}
\left(\int_{\mathbb{R}^{d+k}}|u|^{p^*}dxdy\right)^{\frac{1}{p^*}}\leq C\left(\int_{\mathbb{R}^{d+k}}|\nabla_\mu u|^{p}dxdy\right)^{\frac{1}{p}}.
\end{equation*}
\end{theorem}

Next, we shall prove the associated Hardy-Sobolev type inequalities.
\begin{theorem}\label{Hary-Sobolev}
If $1<p<Q$, $0\leq s\leq p$, there exists a positive constant $C=C(s, p,\mu,Q)$ such that for any $u\in D_0^{1,p}(\mathbb{R}^{d+k})$,
\begin{equation*}
\int_{\mathbb{R}^{d+k}}\frac{|x|^{\mu s}}{\rho^{\mu s}}\frac{|u|^{p_*(s,p,Q)}}{\rho^s}dxdy\leq C\left(\int_{\mathbb{R}^{d+k}}|\nabla_\mu u|^pdxdy\right)^{\frac{Q-s}{Q-p}}.
\end{equation*}
\end{theorem}
In order to prove Theorem \ref{Sobolev-Stein-Embedding}, we consider the fractional integral operator
\begin{equation*}
  I_\nu f(x,y)=\int_{\mathbb{R}^{d+k}}\rho\big((x,y)-(x',y')\big)^{\nu-Q}f(x',y')dx'dy',\quad 0<\nu<Q.
\end{equation*}
The Hardy-Littlewood-Sobolev theorem for $I_\nu$ holds (see \cite{FS}):
\begin{theorem}\label{Hardy-Littlewood-Sobolev} Let $0<\nu<Q$ and $1\leq p<\frac{Q}{\nu}$. Then

(1) if $1<p<Q$, then the condition $\frac{1}{p}-\frac{1}{q}=\frac{\nu}{Q}$ is necessary and sufficient for the boundedness of $I_\nu$ from $L^p(\mathbb{R}^{d+k})$ to $L^q(\mathbb{R}^{d+k})$;

(2) if $p=1$, then the condition $1-\frac{1}{q}=\frac{\nu}{Q}$ is necessary and sufficient for the boundedness of $I_\nu$ from $L^1(\mathbb{R}^{d+k})$ to $L^{q,\infty}(\mathbb{R}^{d+k})$.
\end{theorem}

\noindent\\[2mm]
\noindent{\bbb Proof of Theorem \ref{Sobolev-Stein-Embedding}. } For any $u\in C_0^\infty(\mathbb{R}^{d+k})$, using the integral representation formula for the fundamental solution of $G_\mu$, we have
\begin{equation}\label{representation}
u(x,y)=\int_{\mathbb{R}^{d+k}}\Gamma\big((x,y)-(x',y')\big)G_\mu u(x',y')dx'dy'.
\end{equation}
Noting $G_\mu=\nabla_\mu\cdot\nabla_\mu$ and $\nabla_\mu^*=-\nabla_\mu$ and integrating by parts at the right side of \eqref{representation}, it follows
\begin{equation*}
u(x,y)=\int_{\mathbb{R}^{d+k}}(\nabla_\mu\Gamma)\big((x,y)-(x',y')\big)\nabla_\mu u(x',y')dx'dy'.
\end{equation*}
Since
\begin{equation*}
  |\nabla_\mu\Gamma|=C_\mu|\nabla_\mu\left(\rho^{2-Q}\right)|=C_\mu(Q-2)\rho^{1-Q}|\nabla_\mu\rho|\leq C_\mu(Q-2)\rho^{1-Q},
\end{equation*}
we obtain
\begin{align*}
|u(x,y)|&\leq C_\mu(Q-2)\int_{\mathbb{R}^{d+k}}\rho\big((x,y)-(x',y')\big)^{1-Q}\left|\nabla_\mu u(x',y')\right|dx'dy'\\
&= C_\mu(Q-2)I_1(|\nabla_\mu u|).
\end{align*}
Now, applying Theorem \ref{Hardy-Littlewood-Sobolev}, it yields
\begin{align*}
\left(\int_{\mathbb{R}^{d+k}}|u|^qdxdy\right)^{\frac{1}{q}}
&\leq C_\mu(Q-2)\left(\int_{\mathbb{R}^{d+k}}\left(I_1(|\nabla_\mu u|)\right)^qdxdy\right)^{\frac{1}{q}}\\
&\leq C\left(\int_{\mathbb{R}^{d+k}}|\nabla_\mu u|^{p}dxdy\right)^{\frac{1}{p}},
\end{align*}
where
\begin{equation*}
  \frac{1}{p}-\frac{1}{q}=\frac{1}{Q}\quad \left(\Leftrightarrow q=\frac{Qp}{Q-p}=p^*\right),
\end{equation*}
and $C$ is a suitable positive constant depending only on $p,\mu$ and $Q$. This ends the proof.

\noindent\\[2mm]
\noindent{\bbb Proof of Theorem \ref{Hary-Sobolev}. } Recall $p_*(s,p,Q)=\frac{p(Q-s)}{Q-p}$, where $1<p<Q$ and $0\leq s\leq p$. If $s=0$, then $p_*(p,p,Q)=p$ and Theorem \ref{Hary-Sobolev} is reduced to Theorem \ref{Sobolev-Stein-Embedding}. If $s=p$, then $p_*(0,p,Q)=p^*$ and Theorem \ref{Hary-Sobolev} is reduced to Theorem \ref{Hardy-type-theorem} in the case $\alpha=0$. Theorefore, it suffices to deal with the case $0<s<p$.

Denoting $p_*(s,p,Q)=\left(1-\frac{s}{p}\right)p^*+s$, by H\"{o}lder inquality, \eqref{Hardy-type} in the case $\alpha=0$ and Theorem \ref{Sobolev-Stein-Embedding}, we have
\begin{align*}
\quad &\int_{\mathbb{R}^{d+k}}\frac{|x|^{\mu s}}{\rho^{\mu s}}\frac{|u|^{p_*(s,p,Q)}}{\rho^s}dxdy\\
&\leq\left(\int_{\mathbb{R}^{d+k}}\frac{|x|^{\mu p}}{\rho^{\mu p}}\frac{|u|^p}{\rho^p}dxdy\right)^{\frac{s}{p}} \left(\int_{\mathbb{R}^{d+k}}|u|^{p^*}dxdy\right)^{1-\frac{s}{p}}\\
&\leq \left(\left(\frac{p}{Q-p}\right)^p\int_{\mathbb{R}^{d+k}}|\nabla_\mu u|^pdxdy\right)^{\frac{s}{p}}
\left(C(p,\mu,Q)^{p^*}\left(\int_{\mathbb{R}^{d+k}}|\nabla_\mu u|^pdxdy\right)^{\frac{p^*}{p}}\right)^{1-\frac{s}{p}}\\
&=C\left(\int_{\mathbb{R}^{d+k}}|\nabla_\mu u|^pdxdy\right)^{\frac{Q-s}{Q-p}},
\end{align*}
where
\begin{equation*}
C=\left(\frac{p}{Q-p}\right)^sC(p,\mu,Q)^{\left(1-\frac{s}{p}\right)p^*}.
\end{equation*}
The proof is completed.

\noindent\\[2mm]
To prove Theorem \ref{W-H-S-new}, we shall introduce two results.
\begin{lemma}\label{p-case} (see \cite{JH}) Let $p\geq 1$. For all $\xi_1,\xi_2\in\mathbb{R}^n$, the following inequalities hold:

(1) if $p\leq2$, then
\begin{align*}
|\xi_1+\xi_2|^p-|\xi_1|^p-p|\xi_1|^{p-2}\langle\xi_1,\xi_2\rangle &\leq C(p)|\xi_2|^p,\\
|\xi_2|^p-|\xi_1|^p-p|\xi_1|^{p-2}\langle\xi_1,\xi_2-\xi_1\rangle &\geq C(p)\frac{|\xi_2-\xi_1|^p}{(|\xi_1|+|\xi_2|)^{2-p}};
\end{align*}

(2) if $p>2$, then
\begin{align*}
|\xi_1+\xi_2|^p-|\xi_1|^p-p|\xi_1|^{p-2}\langle\xi_1,\xi_2\rangle &\leq \frac{p(p-1)}{2}(|\xi_1|+|\xi_2|)^{p-2}|\xi_2|^2,\\
|\xi_2|^p-|\xi_1|^p-p|\xi_1|^{p-2}\langle\xi_1,\xi_2-\xi_1\rangle &\geq C(p)|\xi_2-\xi_1|^p,
\end{align*}
where $\langle\cdot,\cdot\rangle$ represents the common inner product on the Euclidean space $\mathbb{R}^n$.
\end{lemma}

\begin{lemma}\label{lambda-case}(see \cite{KJY}) For any $\xi,\eta\in\mathbb{R}^n$ and $\lambda>0$,
\begin{equation*}
|\xi|^{\lambda+1}+\lambda |\eta|^{\lambda+1}-(\lambda+1)|\eta|^{\lambda-1}\langle\xi,\eta\rangle\geq0,
\end{equation*}
and the equality holds if and only if $\xi=\eta$.
\end{lemma}

\noindent\\[2mm]
\noindent{\bbb Proof of Theorem \ref{W-H-S-new}. }  The condition $\alpha>\frac{p-Q}{p}$ implies
\begin{equation*}
\alpha p_*(s,p,Q)-s+Q>0,\quad p\alpha+Q>0,
\end{equation*}
which ensures that the left and right integral of \eqref{weighted-H-S-Grushin} are well defined on $C_0^\infty(\mathbb{R}^{d+k})$.
For any $u\in D_\alpha^{1,p}(\mathbb{R}^{d+k})$, taking $w=\rho^\alpha u$, by the property of convex fuctions and \eqref{gradient-Delta}, we have
\begin{align*}
\left|\nabla_\mu w\right|^p&=\left|\rho^\alpha \nabla_\mu u+\alpha \rho^{\alpha-1} u\nabla_\mu \rho\right|^p\\
                           &\leq\left(\rho^\alpha |\nabla_\mu u|+|\alpha| |\rho^{\alpha-1} u||\nabla_\mu \rho|\right)^p\\
                           &= 2^p\left(\frac{\rho^\alpha |\nabla_\mu u|+|\alpha| |\rho^{\alpha-1} u||\nabla_\mu \rho|}{2}\right)^p\\
                           &\leq 2^p\left(\frac{\rho^{\alpha p} |\nabla_\mu u|^p}{2}+\frac{|\alpha|^p |\rho^{\alpha-1} u|^p|\nabla_\mu \rho|^p}{2}\right)\\
                           &\leq 2^{p-1}\left(\rho^{\alpha p} |\nabla_\mu u|^p+|\alpha|^p \rho^{\alpha p}\frac{|x|^{\mu p}}{\rho^{\mu p}}\frac{|u|^p}{\rho^p}\right).
\end{align*}
It follows from \eqref{Hardy-type} that
\begin{align*}
\int_{\mathbb{R}^{d+k}}\left|\nabla_\mu w\right|^pdxdy&\leq 2^{p-1}\int_{\mathbb{R}^{d+k}}\left(\rho^{\alpha p} |\nabla_\mu u|^p+|\alpha|^p \rho^{\alpha p}\frac{|x|^{\mu p}}{\rho^{\mu p}}\frac{|u|^p}{\rho^p}\right)dxdy\\
                                                      &\leq 2^{p-1}\left(|\alpha|^p\left(\frac{p}{Q-p+\alpha p}\right)^p+1\right)\int_{\mathbb{R}^{d+k}}\rho^{\alpha p} |\nabla_\mu u|^pdxdy,
\end{align*}
which implies $w\in D_0^{1,p}(\mathbb{R}^{d+k})$. In addition, a straightforward computation deduces
\begin{align*}
\int_{\mathbb{R}^{d+k}}\frac{|x|^{\mu s}}{\rho^{\mu s}}\frac{|\rho^\alpha u|^{p_*(s,p,Q)}}{\rho^s}dxdy&=\int_{\mathbb{R}^{d+k}}\frac{|x|^{\mu s}}{\rho^{\mu s}}\frac{|w|^{p_*(s,p,Q)}}{\rho^s}dxdy,\\
\int_{\mathbb{R}^{d+k}}|\rho^\alpha\nabla_\mu u|^pdxdy&=\int_{\mathbb{R}^{d+k}}|\nabla_\mu w-\alpha \rho^{-1} w\nabla_\mu \rho|^pdxdy.
\end{align*}
Therefore, \eqref{weighted-H-S-Grushin} is equivalent to the following inequality
\begin{equation*}
\int_{\mathbb{R}^{d+k}}\frac{|x|^{\mu s}}{\rho^{\mu s}}\frac{|w|^{p_*(s,p,Q)}}{\rho^s}dxdy\leq C\left(\int_{\mathbb{R}^{d+k}}|\nabla_\mu w-\alpha \rho^{-1} w\nabla_\mu \rho|^pdxdy\right)^{\frac{Q-s}{Q-p}}.
\end{equation*}
Noting Theorem \ref{Hary-Sobolev}, it suffices to prove
\begin{equation}\label{Goal-H-S-function}
\int_{\mathbb{R}^{d+k}}|\nabla_\mu w-\alpha \rho^{-1} w\nabla_\mu \rho|^pdxdy\geq C'\int_{\mathbb{R}^{d+k}}|\nabla_\mu w|^pdxdy
\end{equation}
for some suitable constant $C'>0$.\\

According to Lemma \ref{p-case}, we will investigate \eqref{Goal-H-S-function} under the case $1< p\leq 2$ and $p>2$ respectively.

 {\bbb Case 1: }$1< p\leq 2$. Taking $\xi_1=\alpha \rho^{-1} w\nabla_\mu \rho$ and $\xi_2=\nabla_\mu w-\alpha \rho^{-1} w\nabla_\mu \rho$ in the first case of Lemma \ref{p-case}, it indicates
\begin{equation}\label{case-1}
\begin{aligned}
&\quad C(p)\int_{\mathbb{R}^{d+k}}|\nabla_\mu w-\alpha \rho^{-1} w\nabla_\mu \rho|^pdxdy\\
&\leq \int_{\mathbb{R}^{d+k}}\left\{|\nabla_\mu w|^p-|\alpha|^p\frac{|w|^p}{\rho^p}|\nabla_\mu \rho|^p \right.\\
&\qquad\left.-p\alpha |\alpha|^{p-2}\frac{w|w|^{p-2}}{\rho^{p-1}}|\nabla_\mu \rho|^{p-2}\langle\nabla_\mu \rho,\nabla_\mu w-\alpha \rho^{-1} w\nabla_\mu \rho\rangle\right\} dxdy\\
&=\int_{\mathbb{R}^{d+k}}\left\{|\nabla_\mu w|^p+(p-1)|\alpha|^p\frac{|w|^p}{\rho^p}|\nabla_\mu \rho|^p \right.\\
&\qquad\left.-p\alpha |\alpha|^{p-2}\frac{w|w|^{p-2}}{\rho^{p-1}}|\nabla_\mu \rho|^{p-2}\langle\nabla_\mu \rho,\nabla_\mu w\rangle\right\} dxdy\\
&=\int_{\mathbb{R}^{d+k}}\left\{|\nabla_\mu w|^p+(p-1)|\alpha|^p\frac{|x|^{\mu p}}{\rho^{\mu p}}\frac{|w|^p}{\rho^p} \right.\\
&\qquad-\alpha |\alpha|^{p-2}\rho^{1-p}|\nabla_\mu \rho|^{p-2}\langle\nabla_\mu \rho,\nabla_\mu |w|^p\rangle\bigg\} dxdy.
\end{aligned}
\end{equation}

Integrating by parts, we have
\begin{equation}\label{case-1-integration}
\begin{aligned}
&\quad\int_{\mathbb{R}^{d+k}}\rho^{1-p}|\nabla_\mu \rho|^{p-2}\langle\nabla_\mu \rho,\nabla_\mu |w|^p\rangle dxdy\\
&=-\int_{\mathbb{R}^{d+k}}|w|^p \text{div}_\mu\left(\rho^{1-p}|\nabla_\mu \rho|^{p-2}\nabla_\mu \rho\right) dxdy.
\end{aligned}
\end{equation}

Since
\begin{equation*}
 \nabla_\mu(|x|^{(p-2)\mu})\cdot\nabla_\mu\rho=(p-2)\mu\frac{|x|^{p\mu}}{\rho^{2\mu+1}},
\end{equation*}
by \eqref{gradient-Delta}, a straightforward computation implies
\begin{equation}\label{case-1-computation}
\begin{aligned}
 &\quad \text{div}_\mu\left(\rho^{1-p}|\nabla_\mu \rho|^{p-2}\nabla_\mu \rho\right)\\
 &=\text{div}_\mu\left(|x|^{(p-2)\mu}\rho^{1-p+(2-p)\mu}\nabla_\mu \rho\right)\\
 &=\rho^{1-p+(2-p)\mu}\nabla_\mu(|x|^{(p-2)\mu})\cdot\nabla_\mu\rho\\
 &+\big(1-p+(2-p)\mu\big)|x|^{(p-2)\mu}\rho^{(2-p)\mu-p}|\nabla_\mu \rho|^2+|x|^{(p-2)\mu}\rho^{1-p+(2-p)\mu}G_\mu \rho\\
 &=(Q-p)\frac{|x|^{p\mu}}{\rho^{(p+1)\mu}}.
\end{aligned}
\end{equation}

Putting \eqref{case-1-integration} and \eqref{case-1-computation} into \eqref{case-1}, , we have
\begin{equation}\label{last-case-1}
\begin{aligned}
&C(p)\int_{\mathbb{R}^{d+k}}|\nabla_\mu w-\alpha \rho^{-1} w\nabla_\mu \rho|^pdxdy\geq \int_{\mathbb{R}^{d+k}}|\nabla_\mu w|^pdxdy\\
&\qquad \qquad +\alpha |\alpha|^{p-2}\big(Q-p+(p-1)\alpha\big)\int_{\mathbb{R}^{d+k}}\frac{|x|^{\mu p}}{\rho^{\mu p}}\frac{|w|^p}{\rho^p}dxdy.
\end{aligned}
\end{equation}
Note that the condition $\alpha>\frac{p-Q}{p}$ deduces that $Q-p+(p-1)\alpha>0$.

If $\alpha\leq 0$, it follows from \eqref{Hardy-type} and \eqref{last-case-1} that
\begin{equation}\label{final-case-1}
C(p)\int_{\mathbb{R}^{d+k}}|\nabla_\mu w-\alpha \rho^{-1} w\nabla_\mu \rho|^pdxdy\geq C_1\int_{\mathbb{R}^{d+k}}|\nabla_\mu w|^pdxdy,
\end{equation}
where
\begin{align*}
C_1&=1+\alpha |\alpha|^{p-2}\big(Q-p+(p-1)\alpha\big)\left(\frac{p}{Q-p}\right)^p\\
   &=\left(\frac{p}{Q-p}\right)^p\left\{\left(\frac{Q-p}{p}\right)^p+(p-1)|\alpha|^p+(Q-p)\alpha |\alpha|^{p-2}\right\}.
\end{align*}
Taking $\xi=\frac{p-Q}{p}$, $\eta=\alpha$ and $\lambda=p-1>0$ in Lemma \ref{lambda-case}, it shows $C_1>0$.

If $\alpha>0$, then \eqref{final-case-1} holds naturally with $C_1=1$.

In conclusion, we prove \eqref{Goal-H-S-function} with $C'=C(p)^{-1}C_1>0$.

{\bbb Case 2: }$p>2$. A direct calculation gives
\begin{equation*}
\begin{aligned}
&\quad 2^{p-2}\frac{p(p-1)}{2}\big(|\alpha \rho^{-1} w\nabla_\mu \rho|+|\nabla_\mu w|\big)^{p-2}|\nabla_\mu w-\alpha \rho^{-1} w\nabla_\mu \rho|^2\\
&\geq \frac{p(p-1)}{2}\big(2|\alpha \rho^{-1} w\nabla_\mu \rho|+|\nabla_\mu w|\big)^{p-2}|\nabla_\mu w-\alpha \rho^{-1} w\nabla_\mu \rho|^2\\
&\geq \frac{p(p-1)}{2}\big(|\alpha \rho^{-1} w\nabla_\mu \rho|+|\nabla_\mu w-\alpha \rho^{-1} w\nabla_\mu \rho|\big)^{p-2}|\nabla_\mu w-\alpha \rho^{-1} w\nabla_\mu \rho|^2.
\end{aligned}
\end{equation*}
As Case 1, applying the estimate in the second case of Lemma \ref{p-case} to the above inequality deduces
\begin{align*}
&\quad 2^{p-2}\frac{p(p-1)}{2}\big(|\alpha \rho^{-1} w\nabla_\mu \rho|+|\nabla_\mu w|\big)^{p-2}|\nabla_\mu w-\alpha \rho^{-1} w\nabla_\mu \rho|^2\\
&\geq |\nabla_\mu w|^p-|\alpha \rho^{-1} w\nabla_\mu \rho|^p-p|\alpha \rho^{-1} w\nabla_\mu \rho|^{p-2}\langle\alpha \rho^{-1} w\nabla_\mu \rho,\nabla_\mu w-\alpha \rho^{-1} w\nabla_\mu \rho\rangle\\
&=|\nabla_\mu w|^p+(p-1)|\alpha|^p\frac{|x|^{\mu p}}{\rho^{\mu p}}\frac{|w|^p}{\rho^p}-\alpha |\alpha|^{p-2}\rho^{1-p}|\nabla_\mu \rho|^{p-2}\langle\nabla_\mu \rho,\nabla_\mu |w|^p\rangle.
\end{align*}
Therefore, argued as Case 1,
\begin{equation}\label{step-1-case-2}
\begin{aligned}
&2^{p-2}\frac{p(p-1)}{2}\int_{\mathbb{R}^{d+k}}\big(|\alpha \rho^{-1} w\nabla_\mu \rho|+|\nabla_\mu w|\big)^{p-2}|\nabla_\mu w-\alpha \rho^{-1} w\nabla_\mu \rho|^2dxdy\\
&\geq \int_{\mathbb{R}^{d+k}}\left\{|\nabla_\mu w|^p+(p-1)|\alpha|^p\frac{|x|^{\mu p}}{\rho^{\mu p}}\frac{|w|^p}{\rho^p}\right.\\
&\qquad -\alpha |\alpha|^{p-2}\rho^{1-p}|\nabla_\mu \rho|^{p-2}\langle\nabla_\mu \rho,\nabla_\mu |w|^p\rangle\bigg\}dxdy\\
&=\int_{\mathbb{R}^{d+k}}|\nabla_\mu w|^pdxdy+\alpha |\alpha|^{p-2}\big(Q-p+(p-1)\alpha\big)\int_{\mathbb{R}^{d+k}}\frac{|x|^{\mu p}}{\rho^{\mu p}}\frac{|w|^p}{\rho^p}dxdy\\
&\geq C_1\int_{\mathbb{R}^{d+k}}|\nabla_\mu w|^pdxdy,
\end{aligned}
\end{equation}
holds for a sutitable $C_1>0$.

In addition, exploiting H\"{o}lder inequality and Minkowski inequality, it deduces
\begin{equation}\label{step-21-case-2}
\begin{aligned}
&\quad \int_{\mathbb{R}^{d+k}}\big(|\alpha \rho^{-1} w\nabla_\mu \rho|+|\nabla_\mu w|\big)^{p-2}|\nabla_\mu w-\alpha \rho^{-1} w\nabla_\mu \rho|^2dxdy\\
&\leq \left(\int_{\mathbb{R}^{d+k}}\big(|\alpha \rho^{-1} w\nabla_\mu \rho|+|\nabla_\mu w|\big)^pdxdy\right)^\frac{p-2}{p}\\
&\qquad \times \left(\int_{\mathbb{R}^{d+k}}|\nabla_\mu w-\alpha \rho^{-1} w\nabla_\mu \rho|^pdxdy\right)^\frac{2}{p}\\
&\leq \left\{\left(\int_{\mathbb{R}^{d+k}}\big|\alpha \rho^{-1} w\nabla_\mu \rho|^pdxdy\right)^\frac{1}{p}+\left(\int_{\mathbb{R}^{d+k}}|\nabla_\mu w|^pdxdy\right)^\frac{1}{p}\right\}^{p-2}\\
&\qquad \times \left(\int_{\mathbb{R}^{d+k}}|\nabla_\mu w-\alpha \rho^{-1} w\nabla_\mu \rho|^pdxdy\right)^\frac{2}{p}.
\end{aligned}
\end{equation}
We conclude from \eqref{Hardy-type},
\begin{equation}\label{step-22-case-2}
\begin{aligned}
&\quad \int_{\mathbb{R}^{d+k}}\big|\alpha \rho^{-1} w\nabla_\mu \rho|^pdxdy\\
&=|\alpha|^p\int_{\mathbb{R}^{d+k}}\frac{|x|^{\mu p}}{\rho^{\mu p}}\frac{|w|^p}{\rho^p}dxdy\\
&\leq|\alpha|^p\left(\frac{p}{Q-p}\right)^p\int_{\mathbb{R}^{d+k}}|\nabla_\mu w|^pdxdy.
\end{aligned}
\end{equation}
Hence, by \eqref{step-21-case-2} and \eqref{step-22-case-2},
\begin{equation}\label{step-2-case-2}
\begin{aligned}
&\quad \int_{\mathbb{R}^{d+k}}\big(|\alpha \rho^{-1} w\nabla_\mu \rho|+|\nabla_\mu w|\big)^{p-2}|\nabla_\mu w-\alpha \rho^{-1} w\nabla_\mu \rho|^2dxdy\\
&\leq \left(|\alpha|\frac{p}{Q-p}+1\right)^{p-2}\left(\int_{\mathbb{R}^{d+k}}|\nabla_\mu w|^pdxdy\right)^\frac{p-2}{p}\\
&\qquad \times \left(\int_{\mathbb{R}^{d+k}}|\nabla_\mu w-\alpha \rho^{-1} w\nabla_\mu \rho|^pdxdy\right)^\frac{2}{p}.
\end{aligned}
\end{equation}
Combining \eqref{step-1-case-2} and \eqref{step-2-case-2}, it follows
\begin{equation*}
\int_{\mathbb{R}^{d+k}}|\nabla_\mu w-\alpha \rho^{-1} w\nabla_\mu \rho|^pdxdy\geq C'\int_{\mathbb{R}^{d+k}}|\nabla_\mu w|^pdxdy,
\end{equation*}
where $C'=2^{\frac{p(3-p)}{2}}\left[p(p-1)\right]^{-\frac{p}{2}}\left(|\alpha|\frac{p}{Q-p}+1\right)^\frac{p(2-p)}{2}C_1^\frac{p}{2}>0$.

In conclusion, \eqref{Goal-H-S-function} is proved.
\noindent\\[4mm]
\section{Proof of Theorem \ref{CKN-Grushin}}
\subsection{Necessity}

{\bbb \quad Necessity of \eqref{dimensional-balance-new}: } Let $0\not\equiv u\in C_0^\infty(\mathbb{R}^{d+k})$ satisfy \eqref{CKN-inequality-new}.  Then $u_\lambda=u\circ\delta_\lambda (\lambda>0)$ also satisfies \eqref{CKN-inequality-new}. A direct computation shows
\begin{align*}
\int_{\mathbb{R}^{d+k}}\left(\frac{|x|^\mu}{\rho^\mu}\right)^{(\alpha-\gamma)r}\rho^{\gamma r}|u_\lambda|^rdxdy&=\lambda^{-\gamma r-Q}\int_{\mathbb{R}^{d+k}}\left(\frac{|x|^\mu}{\rho^\mu}\right)^{(\alpha-\gamma)r}\rho^{\gamma r}|u|^rdxdy,\\
\int_{\mathbb{R}^{d+k}}\rho^{\alpha p}|\nabla_\mu u_\lambda|^p dxdy&=\lambda^{-(\alpha-1)p-Q}\int_{\mathbb{R}^{d+k}}\rho^{\alpha p}|\nabla_\mu u|^p dxdy,\\
\int_{\mathbb{R}^{d+k}}\left(\frac{|x|^\mu}{\rho^\mu}\right)^{(\alpha-\beta)q}\rho^{\beta q}|u_\lambda|^qdxdy&=\lambda^{-\beta q-Q}\int_{\mathbb{R}^{d+k}}\left(\frac{|x|^\mu}{\rho^\mu}\right)^{(\alpha-\beta)q}\rho^{\beta q}|u|^qdxdy.
\end{align*}
Applying $u_\lambda$ to \eqref{CKN-inequality-new}, it follows
\begin{equation*}
\lambda^{-\gamma-\frac{Q}{r}}\leq \lambda^{a\left(-(\alpha-1)-\frac{Q}{p}\right)+(1-a)\left(-\beta-\frac{Q}{q}\right)},
\end{equation*}
which is true for any $\lambda>0$, so the powers of $\lambda$  on the two sides must be equal, i.e.
\begin{equation*}
-\gamma-\frac{Q}{r}=a\left[-(\alpha-1)-\frac{Q}{p}\right]+(1-a)\left(-\beta-\frac{Q}{q}\right),
\end{equation*}
which is exactly \eqref{dimensional-balance-new}.

{\bbb \quad Necessity of \eqref{index}: } Let $0\not\equiv u\in C_0^\infty(B_1)$ satisfy \eqref{CKN-inequality-new}.  Take $(x_0,y_0)\in\mathbb{R}^{d+k}$, $x_0\neq 0$ and for sufficiently large $\lambda>0$, define
\begin{equation*}
u_\lambda(x,y)=u\left((x,y)-\delta_\lambda(x_0,y_0)\right)=u(x-\lambda x_0,y-\lambda^{1+\mu}y_0).
\end{equation*}

Since
\begin{align*}
&\quad \int_{\mathbb{R}^{d+k}}\left(\frac{|x|^\mu}{\rho^\mu}\right)^{(\alpha-\gamma)r}\rho^{\gamma r}|u_\lambda|^rdxdy\\
&=\int_{B_1}\left(\frac{|x+\lambda x_0|}{\rho(x+\lambda x_0,y+\lambda^{1+\mu}y_0)}\right)^{\mu(\alpha-\gamma)r}\rho(x+\lambda x_0,y+\lambda^{1+\mu}y_0)^{\gamma r}|u|^rdxdy\\
&=\lambda^{\gamma r} \int_{B_1}\left(\frac{|\lambda^{-1}x+x_0|}{\rho(\lambda^{-1}x+x_0,\lambda^{-1}y+y_0)}\right)^{\mu(\alpha-\gamma)r}
\rho(\lambda^{-1}x+x_0,\lambda^{-1}y+y_0)^{\gamma r}|u|^rdxdy\\
&\geq C_1\lambda^{\gamma r},
\end{align*}
\begin{equation*}
\int_{\mathbb{R}^{d+k}}\left(\frac{|x|^\mu}{\rho^\mu}\right)^{(\alpha-\beta)q}\rho^{\beta q}|u_\lambda|^qdxdy\leq C_3\lambda^{\beta q},
\end{equation*}
and
\begin{align*}
&\quad \int_{\mathbb{R}^{d+k}}\rho^{\alpha p}|\nabla_\mu u_\lambda|^p dxdy\\
&=\int_{\mathbb{R}^{d+k}}\rho(x+\lambda x_0,y+\lambda^{1+\mu}y_0)^{\alpha p}|\nabla_\mu u|^p dxdy\\
&=\lambda^{\alpha p}\int_{\mathbb{R}^{d+k}}\rho(\lambda^{-1}x+x_0,\lambda^{-1}y+y_0)^{\alpha p}|\nabla_\mu u|^p dxdy.\\
&\leq C_2\lambda^{\alpha p}
\end{align*}，
applying $u_\lambda$ to \eqref{CKN-inequality-new}, we have
\begin{equation*}
C_1^{\frac{1}{r}}\lambda^\gamma\leq C_2^{\frac{a}{p}}C_3^{\frac{1-a}{q}}\lambda^{a\alpha+(1-a)\beta},
\end{equation*}
which yields $\gamma=a\sigma+(1-a)\beta\leq a\alpha+(1-a)\beta$, namely \eqref{index}.

{\bbb \quad Necessity of \eqref{case>1}: }  We conclude from \eqref{dimensional-balance-new} that
\begin{equation*}
\frac{1}{p}+\frac{\alpha-1}{Q}=\frac{1}{q}+\frac{\beta}{Q}=\frac{1}{r}+\frac{\gamma}{Q}.
\end{equation*}

Choose the function
\begin{equation*}
u_\varepsilon=\left\{
\begin{aligned}
&0\qquad \qquad \qquad \text{for } \rho\geq1\\
&\rho^{-\gamma-\frac{Q}{r}}\log\frac{1}{\rho}\quad \text{for } \varepsilon\leq\rho\leq1\\
&\varepsilon^{-\gamma-\frac{Q}{r}}\log\frac{1}{\varepsilon}\quad \text{for } \rho\leq\varepsilon.
\end{aligned}
\right.
\end{equation*}
By polar coordinate changes, it leads to
\begin{align*}
&\quad \int_{\mathbb{R}^{d+k}}\left(\frac{|x|^\mu}{\rho^\mu}\right)^{(\alpha-\gamma)r}\rho^{\gamma r}|u_\varepsilon|^rdxdy\\
&=C'\left(\int_0^\varepsilon \rho^{\gamma r+Q-1}\left(\varepsilon^{-\gamma-\frac{Q}{r}}\log\frac{1}{\varepsilon}\right)^rd\rho+\int_\varepsilon^1\rho^{-1}\left(\log\frac{1}{\rho}\right)^rd\rho\right)\\
&=C'\left(\frac{\left(\log\frac{1}{\varepsilon}\right)^r}{\gamma r+Q}+\frac{\left(\log\frac{1}{\varepsilon}\right)^{r+1}}{r+1}\right)\leq C_1\left(\log\frac{1}{\varepsilon}\right)^{r+1},
\end{align*}
\begin{equation*}
\int_{\mathbb{R}^{d+k}}\left(\frac{|x|^\mu}{\rho^\mu}\right)^{(\alpha-\beta)q}\rho^{\beta q}|u_\varepsilon|^qdxdy\leq C_2\left(\log\frac{1}{\varepsilon}\right)^{q+1},
\end{equation*}
and
\begin{align*}
&\quad \int_{\mathbb{R}^{d+k}}\rho^{\alpha p}|\nabla_\mu u_\varepsilon|^p dxdy\\
&\leq2^{p-1}\int_{\varepsilon\leq\rho\leq1}\rho^{-Q}\left(\left(\frac{Q-p}{p}\right)^p\log^p\frac{1}{\rho}+1\right)|\nabla_\mu \rho|^p dxdy\\
&\leq C_3\left(\log\frac{1}{\varepsilon}\right)^{p+1}.
\end{align*}
Applying $u_\varepsilon$ to \eqref{CKN-inequality-new}, we have
\begin{equation*}
C_1^\frac{1}{r}\left(\log\frac{1}{\varepsilon}\right)^{1+\frac{1}{r}}\leq C_2^\frac{a}{p}C_3^\frac{1-a}{q} \left(\log\frac{1}{\varepsilon}\right)^{a\left(1+\frac{1}{p}\right)+(1-a)\left(1+\frac{1}{q}\right)},
\end{equation*}
which implies that
\begin{equation*}
1+\frac{1}{r}\leq a\left(1+\frac{1}{p}\right)+(1-a)\left(1+\frac{1}{q}\right),
\end{equation*}
which immediately leads to
\begin{equation}\label{case>1-step}
\frac{1}{r}\leq \frac{a}{p}+\frac{1-a}{q}.
\end{equation}
Combining \eqref{alpha-beta-gamma}, \eqref{dimensional-balance-new} and \eqref{case>1-step} yields \eqref{case>1}.

\noindent\\[4mm]
\subsection{Sufficiency}
If $a=0$, \eqref{CKN-inequality-new} holds true obviously; If $a=1$, the proof is complete in Remark \ref{remark}. In the sequel, we deal only with the case $0<a<1$.

{\bbb Case $\uppercase\expandafter{\romannumeral 1}$: $0<a<1$, $0\leq\alpha-\sigma\leq1$}. In this case, $p\leq\left(\frac{1}{p}+\frac{\alpha-\sigma-1}{Q}\right)^{-1}\leq p^*$. Analogous to the argument in Remark \ref{remark}, there exists $t\in [0,1]$ satisfying
$\left(\frac{1}{p}+\frac{\alpha-\sigma-1}{Q}\right)^{-1}=\frac{p(Q-tp)}{Q-p}$ and $(\alpha-\sigma)\left(\frac{1}{p}+\frac{\alpha-\sigma-1}{Q}\right)^{-1}=tp$. Applying \eqref{weighted-H-S-Grushin} with $0\leq s=tp\leq p$, we obtain
\begin{equation}\label{part-one}
\begin{aligned}
&\quad \int_{\mathbb{R}^{d+k}}\left(\frac{|x|^\mu}{\rho^\mu}\right)^{(\alpha-\sigma)\left(\frac{1}{p}+\frac{\alpha-\sigma-1}{Q}\right)^{-1}}\frac{|\rho^\alpha u|^{\left(\frac{1}{p}+\frac{\alpha-\sigma-1}{Q}\right)^{-1}}}{\rho^{(\alpha-\sigma)\left(\frac{1}{p}+\frac{\alpha-\sigma-1}{Q}\right)^{-1}}}dxdy\\
&\leq C\left(\int_{\mathbb{R}^{d+k}}|\rho^\alpha\nabla_\mu u|^pdxdy\right)^{\frac{Q-tp}{Q-p}}.
\end{aligned}
\end{equation}
From \eqref{alpha-beta-gamma} and \eqref{dimensional-balance-new},
\begin{align*}
\frac{1}{r}&=a\left(\frac{1}{p}+\frac{\alpha-\sigma-1}{Q}\right)+\frac{1-a}{q}<1,\quad(\Rightarrow r>1),\\
\left(\frac{|x|^\mu}{\rho^\mu}\right)^{\alpha-\gamma}\rho^\gamma|u|&=\left\{\left(\frac{|x|^\mu}{\rho^\mu}\right)^{a(\alpha-\sigma)}\rho^{a\sigma}|u|^a\right\}\left\{\left(\frac{|x|^\mu}{\rho^\mu}\right)^{(1-a)(\alpha-\beta)}\rho^{(1-a)\beta}|u|^{1-a}\right\}.
\end{align*}
By H\"{o}lder inequality and \eqref{part-one}, it follows
\begin{align*}
&\quad \left(\int_{\mathbb{R}^{d+k}}\left(\frac{|x|^\mu}{\rho^\mu}\right)^{(\alpha-\gamma)r}\rho^{\gamma r}|u|^rdxdy\right)^{\frac{1}{r}} \\
&\leq \left(\int_{\mathbb{R}^{d+k}}\left(\frac{|x|^\mu}{\rho^\mu}\right)^{(\alpha-\sigma)\left(\frac{1}{p}+\frac{\alpha-\sigma-1}{Q}\right)^{-1}}\frac{|\rho^\alpha u|^{\left(\frac{1}{p}+\frac{\alpha-\sigma-1}{Q}\right)^{-1}}}{\rho^{(\alpha-\sigma)\left(\frac{1}{p}+\frac{\alpha-\sigma-1}{Q}\right)^{-1}}}dxdy\right)^{a\left(\frac{1}{p}+\frac{\alpha-\sigma-1}{Q}\right)}\\
&\qquad \qquad \times \left(\int_{\mathbb{R}^{d+k}}\left(\frac{|x|^\mu}{\rho^\mu}\right)^{(\alpha-\beta)q}\rho^{\beta q}|u|^qdxdy\right)^{\frac{1-a}{q}}\\
&\leq C\left(\int_{\mathbb{R}^{d+k}}\rho^{\alpha p}|\nabla_\mu u|^p dxdy\right)^{\frac{a}{p}}
\left(\int_{\mathbb{R}^{d+k}}\left(\frac{|x|^\mu}{\rho^\mu}\right)^{(\alpha-\beta)q}\rho^{\beta q}|u|^qdxdy\right)^{\frac{1-a}{q}}.
\end{align*}
\\[0.1cm]
\indent{\bbb Case $\uppercase\expandafter{\romannumeral 2}$: $0<a<1$, $\alpha-\sigma>1$}.
Putting
\begin{equation*}
A^p=\int_{\mathbb{R}^{d+k}}|\rho^\alpha\nabla_Lu|^p dxdy,\,\,B^q=\int_{\mathbb{R}^{d+k}}\left(\frac{|x|^\mu}{\rho^\mu}\right)^{(\beta-\alpha)q}\rho^{\beta q}|u|^q dxdy,
\end{equation*}
we see that \eqref{CKN-inequality-new} can be written as
\begin{equation*}
\left(\int_{\mathbb{R}^{d+k}}\left(\frac{|x|^\mu}{\rho^\mu}\right)^{(\alpha-\gamma)r}\rho^{\gamma r}|u|^rdxdy\right)^{\frac{1}{r}}\leq CA^aB^{1-a}.
\end{equation*}
Rescaling $u$ such that $A^aB^{1-a}=1$, our aim becomes to prove
\begin{equation} \label{aim-function}
\left(\int_{\mathbb{R}^{d+k}}\left(\frac{|x|^\mu}{\rho^\mu}\right)^{(\alpha-\gamma)r}\rho^{\gamma r}|u|^rdxdy\right)^{\frac{1}{r}}\leq C.
\end{equation}

Note in Case $\uppercase\expandafter{\romannumeral 1}$, \eqref{CKN-inequality-new} has been proved for $\alpha-\sigma=0$ and $\alpha-\sigma=1$. Therefore,
\begin{align*}
\left(\int_{\mathbb{R}^{d+k}}\left(\frac{|x|^\mu}{\rho^\mu}\right)^{(\alpha-\delta)s}\rho^{\delta s}|u|^sdxdy\right)^{\frac{1}{s}}&\leq C,\\
\left(\int_{\mathbb{R}^{d+k}}\left(\frac{|x|^\mu}{\rho^\mu}\right)^{(\alpha-\epsilon)t}\rho^{\epsilon t}|u|^tdxdy\right)^{\frac{1}{t}}&\leq C,
\end{align*}
provided that $\delta,s,\epsilon$ and $t$ satisfy
\begin{equation*}
\begin{aligned}
\delta&=b\alpha+(1-b)\beta,\\
\frac{1}{s}&=\frac{b}{p}+\frac{1-b}{q}-\frac{b}{Q},\\
\epsilon&=c(\alpha-1)+(1-c)\beta,\\
\frac{1}{t}&=\frac{c}{p}+\frac{1-c}{q},
\end{aligned}
\end{equation*}
for some choices of $b$ and $c$, $0\leq b,c\leq1$ and
\begin{equation*}
\begin{aligned}
d+\mu(\alpha-\delta)s&>0,\,\delta s+Q>0,\\
d+\mu(\alpha-\epsilon)t&>0,\,\epsilon t+Q>0.
\end{aligned}
\end{equation*}

One computes that
\begin{align*}
\frac{1}{t}+\frac{\epsilon}{Q}&=c\left(\frac{1}{p}+\frac{\alpha-1}{Q}\right)+(1-c)\left(\frac{1}{q}+\frac{\beta}{Q}\right),\\
\frac{1}{r}+\frac{\gamma}{Q}&=a\left(\frac{1}{p}+\frac{\alpha-1}{Q}\right)+(1-a)\left(\frac{1}{q}+\frac{\beta}{Q}\right),\\
\frac{1}{s}+\frac{\delta}{Q}&=b\left(\frac{1}{p}+\frac{\alpha-1}{Q}\right)+(1-b)\left(\frac{1}{q}+\frac{\beta}{Q}\right).
\end{align*}

Since $0<a<1$, $\alpha-\sigma>1$, it ensures $\frac{1}{p}+\frac{\alpha-1}{Q}\not=\frac{1}{r}+\frac{\gamma}{Q}$ from \eqref{case>1} and then $\frac{1}{p}+\frac{\alpha-1}{Q}\not=\frac{1}{q}+\frac{\beta}{Q}$.

1): $\frac{1}{p}+\frac{\alpha-1}{Q}<\frac{1}{q}+\frac{\beta}{Q}$. Taking $b<a<c$, we have
\begin{equation}\label{Q-inquality}
\frac{1}{t}+\frac{\epsilon}{Q}<\frac{1}{r}+\frac{\gamma}{Q}<\frac{1}{s}+\frac{\delta}{Q}.
\end{equation}
A direct computation shows
\begin{align*}
\frac{1}{r}-\frac{1}{s}&=(a-b)\left(\frac{1}{p}-\frac{1}{q}-\frac{1}{Q}\right)+\frac{a(\alpha-\sigma)}{Q},\\
\frac{1}{r}-\frac{1}{t}&=(a-c)\left(\frac{1}{p}-\frac{1}{q}\right)+\frac{a(\alpha-\sigma-1)}{Q},\\
\end{align*}
and
\begin{align*}
&\quad\left(\frac{1}{r}+\frac{\mu(\alpha-\gamma)}{d}\right)-\left(\frac{1}{t}+\frac{\mu(\alpha-\epsilon)}{d}\right)\\
&=(a-c)\left(\frac{1}{p}-\frac{1}{q}+\frac{\mu}{d}(\beta+1-\alpha)\right)+a(\alpha-\sigma-1)\left(\frac{1}{Q}+\frac{\mu}{d}\right),\\
&\quad\left(\frac{1}{r}+\frac{\mu(\alpha-\gamma)}{d}\right)-\left(\frac{1}{s}+\frac{\mu(\alpha-\sigma)}{d}\right)\\
&=(a-b)\left(\frac{1}{p}-\frac{1}{q}-\frac{1}{Q}+\frac{\mu}{d}(\beta-\alpha)\right)+a(\alpha-\sigma)\left(\frac{1}{Q}+\frac{\mu}{d}\right).
\end{align*}
The condition $0<a<1$ and $\alpha-\sigma>1$ imply
\begin{align*}
0<&\frac{a(\alpha-\sigma-1)}{Q}<\frac{a(\alpha-\sigma)}{Q},\\
0<&a(\alpha-\sigma-1)\left(\frac{1}{Q}+\frac{\mu}{d}\right)<a(\alpha-\sigma)\left(\frac{1}{Q}+\frac{\mu}{d}\right),
\end{align*}
and for sufficiently small $|b-a|$ and $|a-c|$,
\begin{equation}\label{d-inequality}
\begin{aligned}
\frac{1}{r}&>\frac{1}{s},\quad \frac{1}{r}+\frac{\mu(\alpha-\gamma)}{d}>\frac{1}{s}+\frac{\mu(\alpha-\sigma)}{d},\\
\frac{1}{r}&>\frac{1}{t},\quad \frac{1}{r}+\frac{\mu(\alpha-\gamma)}{d}>\frac{1}{t}+\frac{\mu(\alpha-\epsilon)}{d}.
\end{aligned}
\end{equation}
Combining \eqref{Q-inquality} and \eqref{d-inequality}, we have
\begin{equation}\label{integration-condition}
\begin{aligned}
d+\frac{\mu(\epsilon-\gamma)tr}{t-r}>0&,\quad \frac{(\gamma-\epsilon)tr}{t-r}+Q>0,\\
d+\frac{\mu(\delta-\gamma)sr}{s-r}>0&,\quad \frac{(\gamma-\delta)sr}{s-r}+Q<0.
\end{aligned}
\end{equation}

Choose a fixed $C_0^\infty(\mathbb{R}^{d+k})$ function $\Phi(x,y)$ $(0\leq\Phi\leq1)$ such that
\begin{equation*}
\Phi(x,y)=\left\{
\begin{aligned}
&1,\,\text{if }\rho(x,y)<1,\\
&0,\,\text{if }\rho(x,y)>2.
\end{aligned}
\right.
\end{equation*}
We shall investigate the left side of \eqref{aim-function} by spliting it to two parts. It obtains by H\"{o}lder inequality that
\begin{equation}\label{zero}
\begin{aligned}
 &\quad \left(\int_{\mathbb{R}^{d+k}}\left(\frac{|x|^\mu}{\rho^\mu}\right)^{(\alpha-\gamma)r}\rho^{\gamma r}\Phi|u|^rdxdy\right)^{\frac{1}{r}}\\
 &\leq\left(\int_{\mathbb{R}^{d+k}}\left(\frac{|x|^\mu}{\rho^\mu}\right)^{(\alpha-\epsilon)t}\rho^{\epsilon t}|u|^tdxdy\right)^{\frac{1}{t}}\\
 &\quad\times\left(\int_{\mathbb{R}^{d+k}}\left(\frac{|x|^\mu}{\rho^\mu}\right)^{\frac{(\epsilon-\gamma)rt}{t-r}}\rho^{\frac{(\gamma-\epsilon)rt}{t-r}}
 \Phi^{\frac{t}{t-r}}dxdy\right)^{\frac{1}{r}-\frac{1}{t}}\\
 &\leq C\left(\int_{B_2}|x|^{\frac{\mu(\epsilon-\gamma)rt}{t-r}}\rho^{(\mu+1)(\gamma-\epsilon)\frac{rt}{t-r}}dxdy\right)^{\frac{1}{r}-\frac{1}{t}}
\end{aligned}
\end{equation}
and
\begin{equation}\label{infty}
\begin{aligned}
 &\quad \left(\int_{\mathbb{R}^{d+k}}\left(\frac{|x|^\mu}{\rho^\mu}\right)^{(\alpha-\gamma)r}\rho^{\gamma r}(1-\Phi)|u|^rdxdy\right)^{\frac{1}{r}}\\
 &\leq\left(\int_{\mathbb{R}^{d+k}}\left(\frac{|x|^\mu}{\rho^\mu}\right)^{(\alpha-\delta)s}\rho^{\delta s}|u|^sdxdy\right)^{\frac{1}{s}}\\
 &\quad\times\left(\int_{\mathbb{R}^{d+k}}\left(\frac{|x|^\mu}{\rho^\mu}\right)^{\frac{(\delta-\gamma)rt}{t-r}}\rho^{\frac{(\gamma-\delta)st}{s-r}}
 (1-\Phi)^{\frac{s}{s-r}}dxdy\right)^{\frac{1}{r}-\frac{1}{s}}\\
 &\leq C\left(\int_{\mathbb{R}^{d+k}\setminus B_1}|x|^{\frac{\mu(\delta-\gamma)rs}{s-r}}\rho^{(\mu+1)(\gamma-\delta)\frac{rs}{s-r}}dxdy\right)^{\frac{1}{r}-\frac{1}{s}}.
\end{aligned}
\end{equation}
Moreover, \eqref{integration-condition} ensures the integrals on the right side in \eqref{zero} and \eqref{infty} are bounded, which easily leads to \eqref{aim-function}.

2): $\frac{1}{p}+\frac{\alpha-1}{Q}>\frac{1}{q}+\frac{\beta}{Q}$. Take $c<a<b$ such that $|c-a|$ and $|a-b|$ are sufficiently small. Now \eqref{Q-inquality}-\eqref{infty} still hold true and then the desired result \eqref{aim-function} is derived.

\noindent\\[4mm]
\noindent\bf{\footnotesize Acknowledgements}\quad\rm
{\footnotesize The authors thank deeply Professor Pengcheng Niu for his suggestion
and encouragement on the paper. They were supported by the National Natural
Science Foundation of China (Grant No. 11371036) and the Fundamental Research Funds for the Central Universities (Grant No. 3102015ZY068).}\\[4mm]

\noindent{\bbb{References}}

\end{document}